\documentclass[10pt, reqno, A4paper]{amsart}

\usepackage[letterpaper,margin=0.75in]{geometry}
\usepackage{amsmath,amssymb,amsthm, epsfig}
\usepackage{mathalfa}
\usepackage{amsfonts}
\usepackage{amsmath}
\usepackage{amssymb}
\usepackage[mathscr]{euscript}

\usepackage{graphicx}
\graphicspath{ {images/} }

\usepackage{wrapfig}

\usepackage{color}

\def \div {\mathrm{div}}

\def \suchthat {\ \big | \ }

\def\dfrac{\displaystyle\frac}
 
\def\epsilon{\varepsilon}

\newtheorem{problem}{Problem}
\theoremstyle{definition}

\theoremstyle{remark}

\numberwithin{equation}{section}

\newcommand{\intav}[1]{\mathchoice {\mathop{\vrule width 6pt height 3 pt depth  -2.5pt
\kern -8pt \intop}\nolimits_{\kern -6pt#1}} {\mathop{\vrule width
5pt height 3  pt depth -2.6pt \kern -6pt \intop}\nolimits_{#1}}
{\mathop{\vrule width 5pt height 3 pt depth -2.6pt \kern -6pt
\intop}\nolimits_{#1}} {\mathop{\vrule width 5pt height 3 pt depth
-2.6pt \kern -6pt \intop}\nolimits_{#1}}}

\title[]{Nonlinear diffusion processes: geometric ideas and beyond}

\author[E.V. Teixeira]{Eduardo V. Teixeira}
\address{University of Central Florida, 4393 Andromeda Loop N, Orlando, FL 32816}{}
\email{eduardo.teixeira@ucf.edu}

\begin{document}

\maketitle
\vspace{-0.6cm}

\begin{center}
{ Professor of Mathematics \\ University of Central Florida \\ Orlando -- Florida, USA.}
\end{center}

\vspace{-2.6cm} 

\hfill  \includegraphics[scale=0.14]{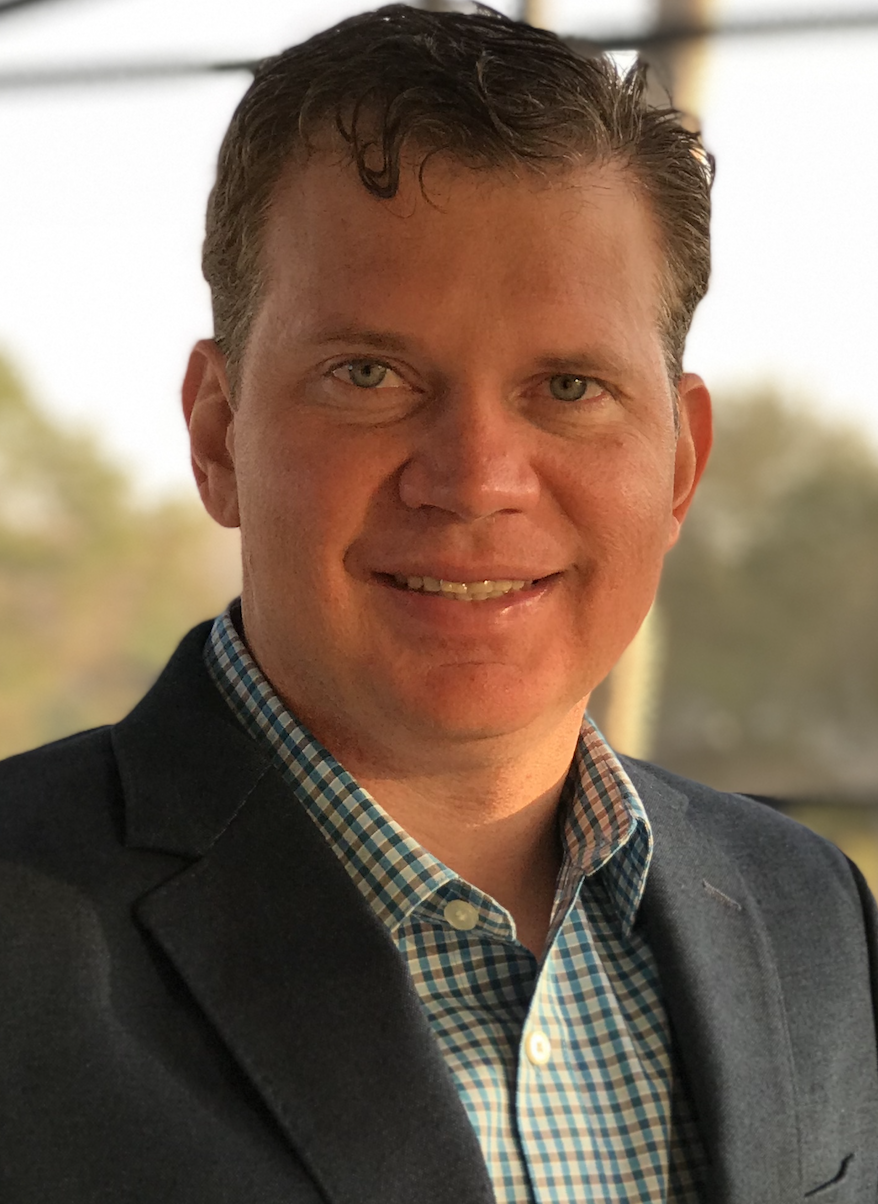}
     
\section{Diffusion}

Diffusion is a phenomenon accounting for average, spread, or balance of quantities in a given process. These constitute innate trends in several fields of natural sciences, which in turn justify why diffusion is such a popular concept among scientists across disciplines. In the realm of mathematics, the study of diffusion is often related to second order differential operators of parabolic type --- or else their stationary versions, the so called elliptic operators.  

The simplest way to appreciate the connection between diffusion and second order elliptic operators is by the following na\"ive looking question: in a domain $\Omega$ of  $\mathbb{R}^n$, find a function $f \colon \Omega \to \mathbb{R}$ such that at each point $y \in \Omega$, $f(y)$ equals the its own average over any ball centered at $y$.  In slightly more precise mathematical terms, we seek for the relation
\begin{equation}\label{avg}
	f(y) = \intav{B_r(y)}{f(x) dx},  
\end{equation}
for all $y \in \Omega$ and all $0<r$ such that $B_r(y) \subset \Omega$. 
The answer may sound surprising at first sight: a function $f$ verifies the averaging property \eqref{avg} if, and only if, it satisfies the so called Laplace equation
\begin{equation}\label{Lap}
	\Delta f(x) : = \frac{\partial^2f(x)}{\partial x_1^2}  +  \frac{\partial^2f(x)}{\partial x_2^2}  + \cdots  \frac{\partial^2f(x)}{\partial x_n^2}  = 0, \quad \forall x \in \Omega.
\end{equation}
The operator appearing in \eqref{Lap} is called the Laplacian, and it is the prototypical example of a second order elliptic operator. Intriguing mysteries  surround the equivalence between \eqref{avg} and \eqref{Lap}. For starters, while \eqref{avg} requires just local integrability of $f$ to make perfect mathematical sense, equation \eqref{Lap} involves second order derivatives of $f$, which in principle, have no reason to exist. This is a key point I want to  emphasize for now; somehow the averaging property \eqref{avg} bears a regularizing effect to a function that verifies it. Understanding this principle is paramount to many different areas of pure and applied mathematics.  

That averages in \eqref{avg} are taken over perfectly symmetric balls conveys the idea of homogeneity of the medium, i.e. there is no preferred direction for diffusion. As for the partial differential equation (PDE) counterpart, \eqref{Lap}, homogeneity translates into an ideal, rotational invariant, constant coefficient operator:  the Laplacian.  

\section{A million ways to say Laplacian} 

Analogy is one of the most powerful features of mathematics and the mathematical theory of diffusion is blessed with many  such correlations, in which problems coming from rather different disciplines lead to a common, unified mathematical treatment. Here is a small pool of samples:
 
\begin{problem} \label{Prob1} What is the equilibrium position of an elastic membrane attached to a given wire?
\end{problem}

Basic physical principles pertaining to the theory of elastic membranes predict that the membrane will adjust itself as to minimize the surface tension. Thus, a first order approximation yields the following minimization problem for the membrane position:
$$
	\min \left \{ \int_{\Omega} |\nabla v|^2 dx \suchthat v = \varphi \text{ on } \partial \Omega \right \}.
$$
The boundary condition $v = \varphi \text{ on } \partial \Omega$ is a way to declare that one is looking for membrane configurations attached to the given wire. It is simple to see that a minimizer of the above functional will satisfy (first in a weak sense and later in the classical sense) the Laplace equation: $\Delta u = 0.$

\begin{problem}\label{Prob2} What is the terminal temperature distribution in a room, prescribed a fixed-in-time wall temperature?
\end{problem}

Let $u(x)$ denote the temperature distribution in the room $\Omega$, prescribed the wall temperature $\varphi \colon \partial \Omega \to \mathbb{R}$. The laws of thermodynamics postulate that the heat flow, $\vec{F}$,  streams from the regions with high  temperature to regions with low temperature. Thus, $\vec{F}$ should be proportional to $-\nabla u$. Let $V\subset \Omega$ be fixed. Since in $V$ no heat is been added nor subtracted, one should have:
$$
	0= \int_{\partial V} \vec{F} \cdot \nu dS = \int_V \text{div} \vec{F} dx = \int_V -\text{div}\left ( \nabla u \right ) dx = \int_V - \Delta u dx.
$$
As $V$ was taken arbitrarily, one finally deduces that the temperature distribution, $u(x)$, in the room $\Omega$, must satisfy:
$\Delta u = 0.$

\begin{problem}\label{Prob3} What is the probability of an 	ant leaving a room through a door before hitting the wall?
\end{problem}

Let $\Omega$ denote the room, $D \subset \partial \Omega$ denote the door and $\varphi \colon \partial \Omega \to \mathbb{R}$ be given by $\varphi(x) = \chi_{D}(x)$, i.e. $\varphi(x) = 1$ if $x \in D$ and $\varphi(x) = 0$ if $x \not \in D$. Let $x \in \Omega$ be the position of the ant and $\delta>0$ the incremental step of the ant towards four possible directions: upwards, downwards, left or right. If $u(x)$ denotes the probability of the ant arriving at the door $D$ before hitting the wall $\partial \Omega \setminus D$, starting from $x \in \Omega$, one can write
$$
	u(x) = \dfrac{u(x + \delta e_1) + u(x - \delta e_1) + u(x + \delta e_2) + u(x - \delta e_2)}{4},
$$
where $e_1 =(1,0)$ and $e_2 = (0,1)$. This is because, being at $x$, the probability of the ant to move from either $x + \delta e_1$ or $x - \delta e_1$ or $x + \delta e_2$ or $x - \delta e_2$ is precisely $\frac{1}{4}$. Dividing the above expression by $\delta^2$ and reorganizing the factors, one reaches:
$$
	0 =  \dfrac{u(x + \delta e_1) + u(x - \delta e_1) - 2u(x)}{\delta^2} + \dfrac{u(x + \delta e_2) + u(x - \delta e_2) - 2u(x)}{\delta^2}.
$$
Letting the incremental step size $\delta$ go to zero, one finds out that the probability $u$ is ruled by the following PDE:
$\Delta u = 0.$

\begin{wrapfigure}{r}{0.5\textwidth}
 \vspace{-1.15cm}
    \includegraphics[scale=0.25]{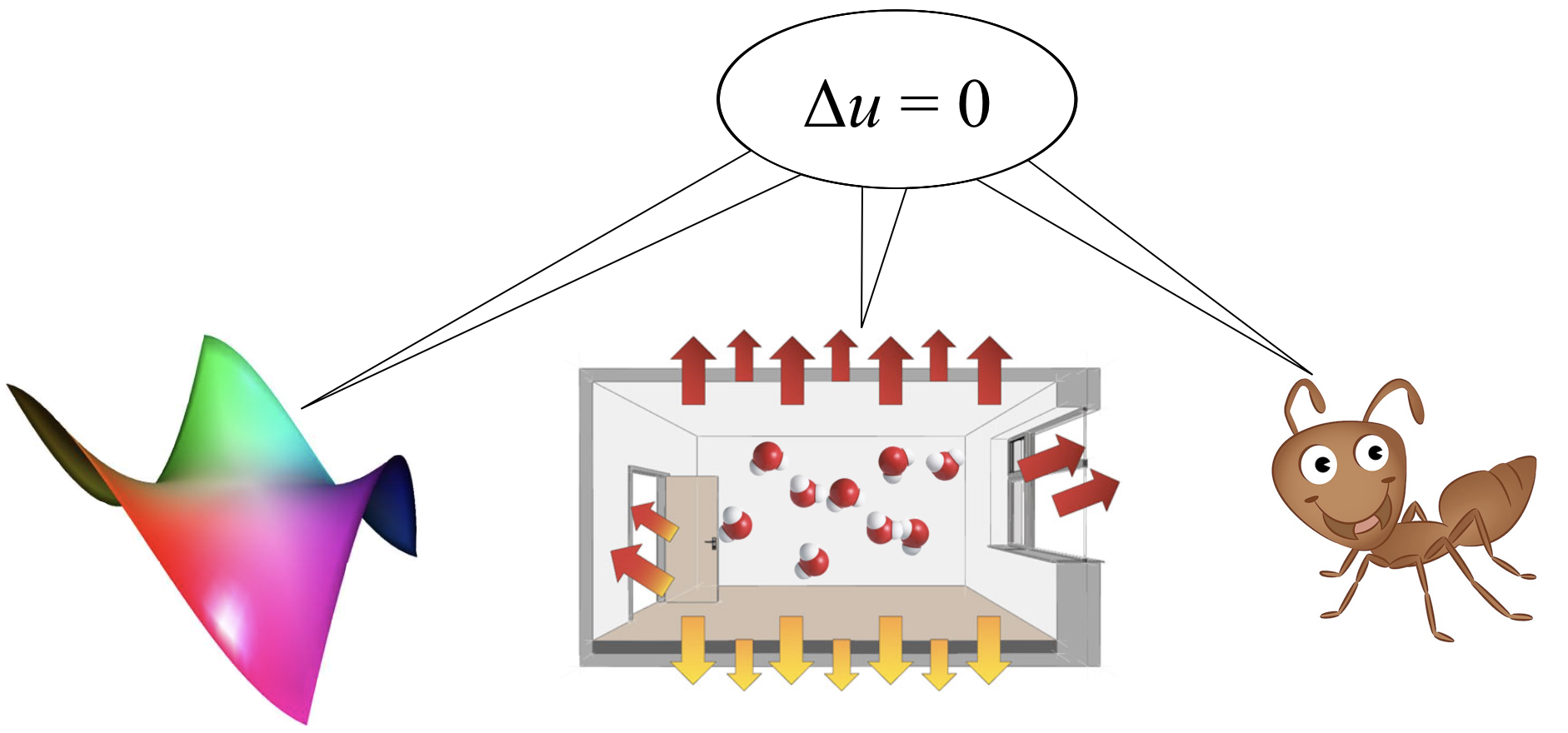}
  \caption{Problems coming from very different realms are linked up through a unified mathematical theory.}
\end{wrapfigure}

\bigskip
 
 \bigskip
 
While I must disclose that rigorous justifications of the above deductions are a bit more laborious, I hope to convey that indeed problems coming from very different backgrounds admit a unified mathematical treatment through the study of the Laplace equation. Even more importantly, such a consolidation yields a bridge between different disciplines, allowing meaningful insight exchanges, which often promote decisive advances in a  field that would hardly be even conjectured otherwise.  \\


\section{Diffusion in complex materials}

More realistic models require more involved differential operators, which may have divergence or non-divergence structures, depending on the nature of the model. Energy considerations, such as in optimization problems or in thermodynamics, often give raise to differential operators in divergence form, whereas probabilistic interpretations of diffusion lead to operators in non-divergence form.  From the mathematical perspective, leading (second order) coefficients  convey the  tangible properties of the medium in which phenomena take place, which, in turn, represent their physical complexities. For instance, in the heat conduction Problem \ref{Prob2}, if one takes into account heterogeneity of the medium ends up with a divergence form equation with non-constant coefficients, say $\div (\gamma(x) \nabla u) = 0$, where $0< \gamma(x) < \infty$.

In some models, not only diffusion is anisotropic, but it can also degenerate at some (a priori unknown) subregions of the domain. Such considerations lead to nonlinear differential operators of degenerate or singular type.  For instance, in the membrane Problem \ref{Prob1}, if instead of minimizing $\int  |\nabla v|^2 dx$, one considers high powers, say $\int  |\nabla v|^p dx$, with $p>2$, then a minimizer will satisfy the so called $p$-Laplace equation, $\text{div} (|\nabla u|^{p-2} \nabla u ) =0$. One should note that, not only the $p$-Laplace operator is non-linear, but actually its ``coefficients", $|\nabla u|^{p-2}$, degenerate along the set of critical points of the solution, $\mathscr{S}:=\{ \nabla u(x) = 0\}$.

Generalizations of the ant Problem \ref{Prob3} yield a mathematical treatment of problems in control and game theory. Emerging equations are of non-divergence form and often involve fully nonlinear structures, $F(x, D^2u)$, that is, are nonlinear with respect to the Hessian argument. In this theory, mathematical manifestation of the diffusion atributes of the operator convert into a monotonicity condition on $F$ with respect to the natural order in the space of symmetric matrices. 

Current literature on the general theory of second order elliptic and parabolic differential operators is vast, dense, and regarded as rather challenging, specially when it comes to understanding regularizing effects of diffusion.  Indeed, analytic approaches to regularity theory mostly involve intricate estimates which are, in general, hard to grasp. As stated by  Mingione (2006), ``Regularity methods are sometimes not very intuitive, and often overburdened by a lot of technical complications, eventually covering the main, basic ideas."

\section{A geometric idea of diffusion}

In contrast to the overwhelming complexity of usual mathematical treatment of nonlinear elliptic operators, if one goes back to the very essence of the idea of diffusion, namely averaging, it becomes  more intuitive that a unified regularity theory could emerge from genuine geometric insights. Thus, (very) loosely speaking, an operator should be considered elliptic if

\begin{center} 

{\it`` it prescribes a balance on how much a solution bends towards each direction."}

\end{center}

Of course, this is not intended to be a mathematical definition {\it per se}; nonetheless, it bears very powerful insights, which, remarkably enough, yield the development of a robust regularity theory solely based upon such a very weak, intuitive notion of averaging process.

\begin{wrapfigure}{r}{0.5\textwidth}
 \begin{center}
 \vspace{-0.4cm}
    \includegraphics[scale=0.22]{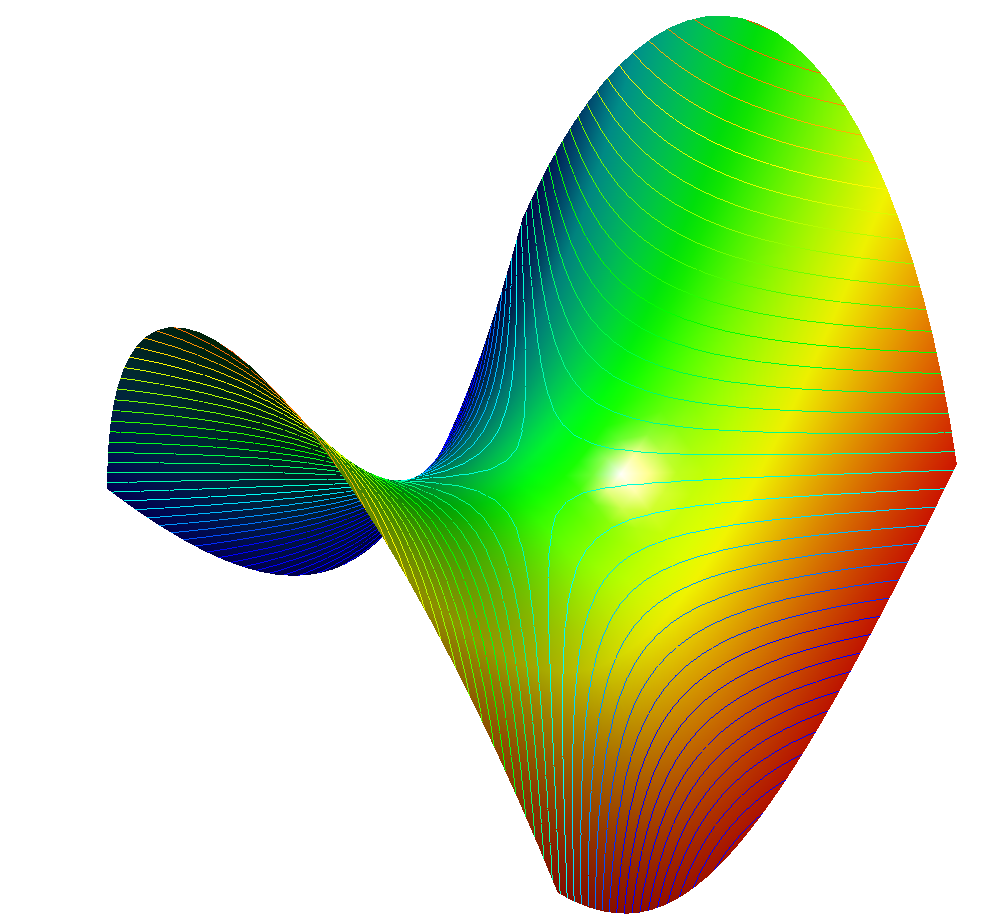}
   \caption{Geometrically speaking, the above function belongs to the class of functions entitled to be a solution of an elliptic equation, as it presents a ``fair" bending balance.}
  \end{center} 
\end{wrapfigure}

The roots of such a radical approach probably go back to Ennio De Giorgi and his magnificent solution to Hilbert's 19th problem, \cite{DeGiorgi}. In a commemorative article, \cite{Bomb},  Enrico Bombieri (1997) mentions a chat he had with De Giorgi on how he got the idea to solve Hilbert's 19th problem. De Giorgi replied as if it was all an indirect consequence of another problem, much more difficult, that he was studying at that moment,  namely the isoperimetric problem in several dimensions. Bombieri records that in his explanation, he kept moving his hands as if he was touching an invisible surface, and showing how to perform his operations and transformations, cutting and pasting invisible masses from one side to the other, leveling and filling the peaks and valleys of theses surfaces. ``I then realized that De Giorgi looked at these functions of several variables literally as geometric objects in space. ... To me, it was an usual way of doing analysis, a field that often requires the use of rather fine estimates, that the normal mathematician grasps more easily through the formulas than through the geometry", comments Bombieri in \cite{Bomb}. He concludes by saying  ``... Perhaps the only other mathematician I met with a geometrical intuition similar to that of De Giorgi's was Luis Caffarelli, of whom De Giorgi was a friend and had a deep esteem." I will come back to Caffarelli in a moment.

While the solution to Hilbert's 19th problem was a major event, it  turned out to be a mere manifestation of a much greater intellectual endowment;  the foundation of De Giorgi's theory of minimal surfaces. This constitutes a rather successful theory developed by De Giorgi and collaborators in the 1960's, where a weak notion of perimeter yields a geometric-measure treatment of the classical Plateau problem, of minimizing area given a prescribed  boundary. It is a parallel endeavor to the famous, and equality successful, Federer--Fleming program, launched in \cite{FF}.

\section{Flatness implies regularity}

  One of the supporting  pillars of De Giorgi's theory of minimal surfaces is the method of {\it flatness improvement}, \cite{DeGiorgiFlat}, which states that if a minimal surface $S$ is ``flat enough", say in $B_1$, with respect to a direction $\nu$, then in $B_{1/2}$, $S$ is even flatter, probably with respect to a slightly tilted direction $\nu'$. Heuristically, the proof of such result goes as follows: suppose, seeking a contradiction, the result is not true. That is, there exists a sequence of minimal surfaces $S_j$ in $B_1$, that are $1/j$ flat with respect to a direction $\nu_j$; however the prospective flatness improvement is not verified in $B_{1/2}$. By {\it compactness}, an appropriate scaling of $S_j$ converges to the graph of a function $f$. By the minimality of $S_j$, $f$ turns out to be a harmonic function, i.e., $\Delta f  = 0$. Being very smooth, the limiting function $f$ does verify the flatness hypothesis. Hence, for $j_0$ sufficiently large, one reaches a contraction on the assumption that no flatness improvement were possible for $S_{j_0}$.


\begin{wrapfigure}{r}{0.5\textwidth}
 \begin{center}
  \vspace{-0.4cm}
    \includegraphics[scale=0.25]{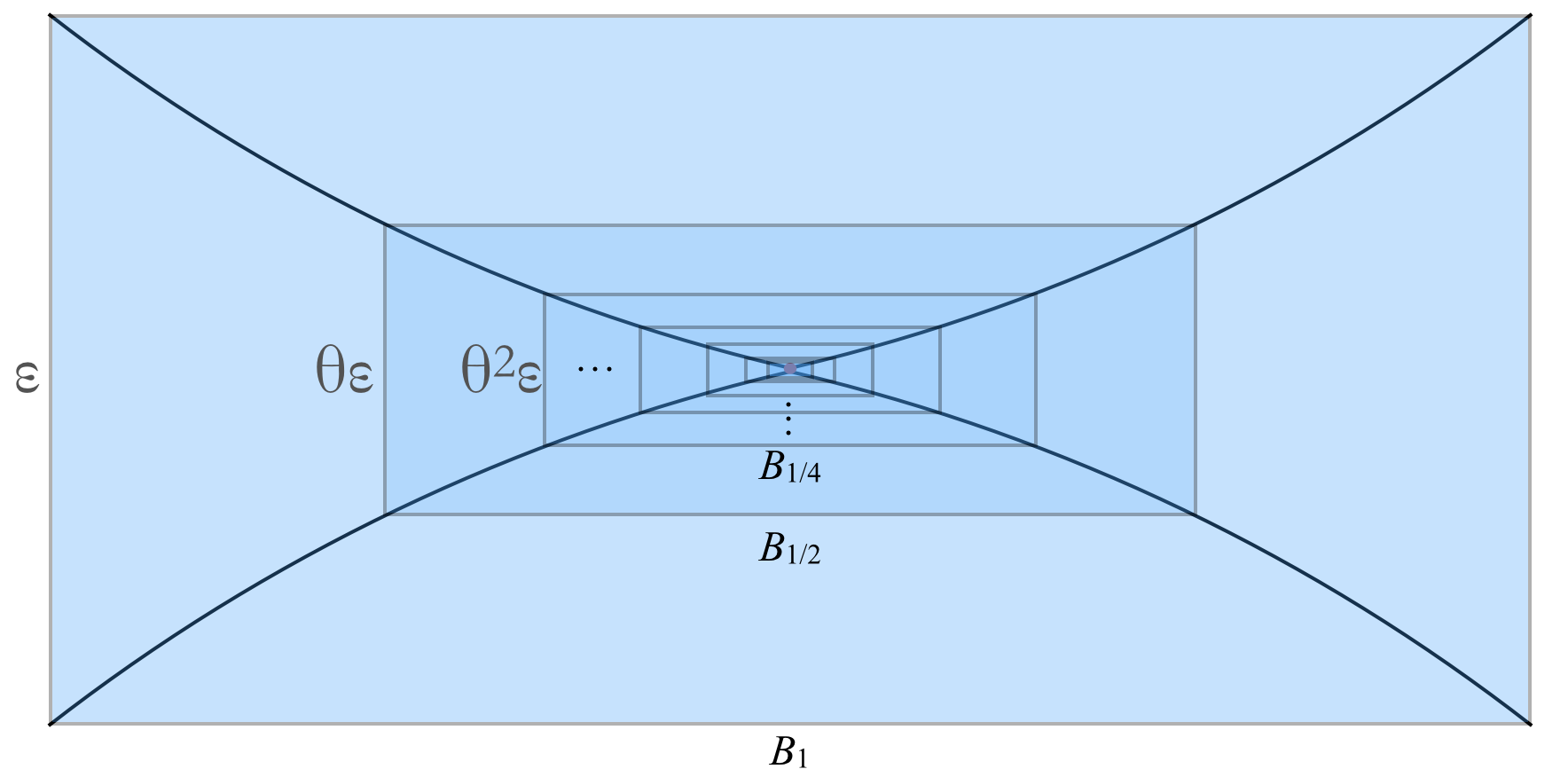}
   \caption{Flatness improvement yields regularity: as $S$ enters $B_{1/2^j}$ it is trapped within a strip of width $\theta^{j} \cdot \varepsilon$, for universal numbers $\varepsilon$ and $\theta$. Ultimately this yields a pathway that conducts $S$ to pass through the origin is a $C^{1,\alpha}$--smooth fashion.}
  \end{center} 
\end{wrapfigure}
The flatness improvement result explicated above is instrumental to ultimately prove that {\it flat enough minimal surfaces are smooth}.   Indeed, in \cite{DeGiorgiFlat}, De Giorgi shows that if a minimal surface $S$ is ``flat enough", say in $B_1$, then in $B_{1/2}$ it is the graph of a $C^{1,\alpha}$ function.

The {motto} then becomes {\it flatness implies regularity}. Little did we know that such a slogan would propagate in many different branches of mathematical analysis, the theory of free boundary problems being one of them.  \\

\section{Free boundaries}

De Giorgi's core ideas and geometric insights were particularly important to the development of the variational theory of free boundary problems. Free boundaries are mathematical manifestations of sharp changes in the parameters that describe the problem. Typically, different physical laws are to be prescribed in distinct, a priori unknown subregions of a domain. This is the case, for instance, of problems involving interfaces between materials, different states of matter, etc.  Free boundaries also arise in physical reactions where interfaces retain some portion of the system's energy, viz., latent heat, membranes, dead cores, flux balances, and so forth. Thus, mathematical models of free boundary problems typically require weak formulations as to give notion to differentiable operations defined on a priori merely measurable sets, and  hence De Giorigi's geometric measure theory is a perfect fit for such endeavor. 

Following up pioneering works of H. Lewy, G. Stampacchia, J.L. Lions, D. Kinderlehrer, among other eminent mathematicians, Luis Caffarelli was the leading figure to pursue a systematic geometric approach to investigate free boundary problems. Caffarelli's 1977 article, \cite{Caff77},  on free boundary regularity for the obstacle problem is a trademark in the theory.  The problem asks for the equilibrium position of an elastic membrane, $u$, restricted to lay above a given obstacle $\psi(x)$. That is, the obstacle problem is the membrane problem \ref{Prob1} with the extra condition of $u$ laying above the obstacle:
$$
	\min \left \{ \int_{\Omega} |\nabla v|^2 dx \suchthat v = \varphi \text{ on } \partial \Omega \text{ and } u(x) \ge \psi(x) \text{ in } \Omega \right \}.
$$
\begin{wrapfigure}{r}{0.5\textwidth}
 \begin{center}
    \includegraphics[scale=0.3]{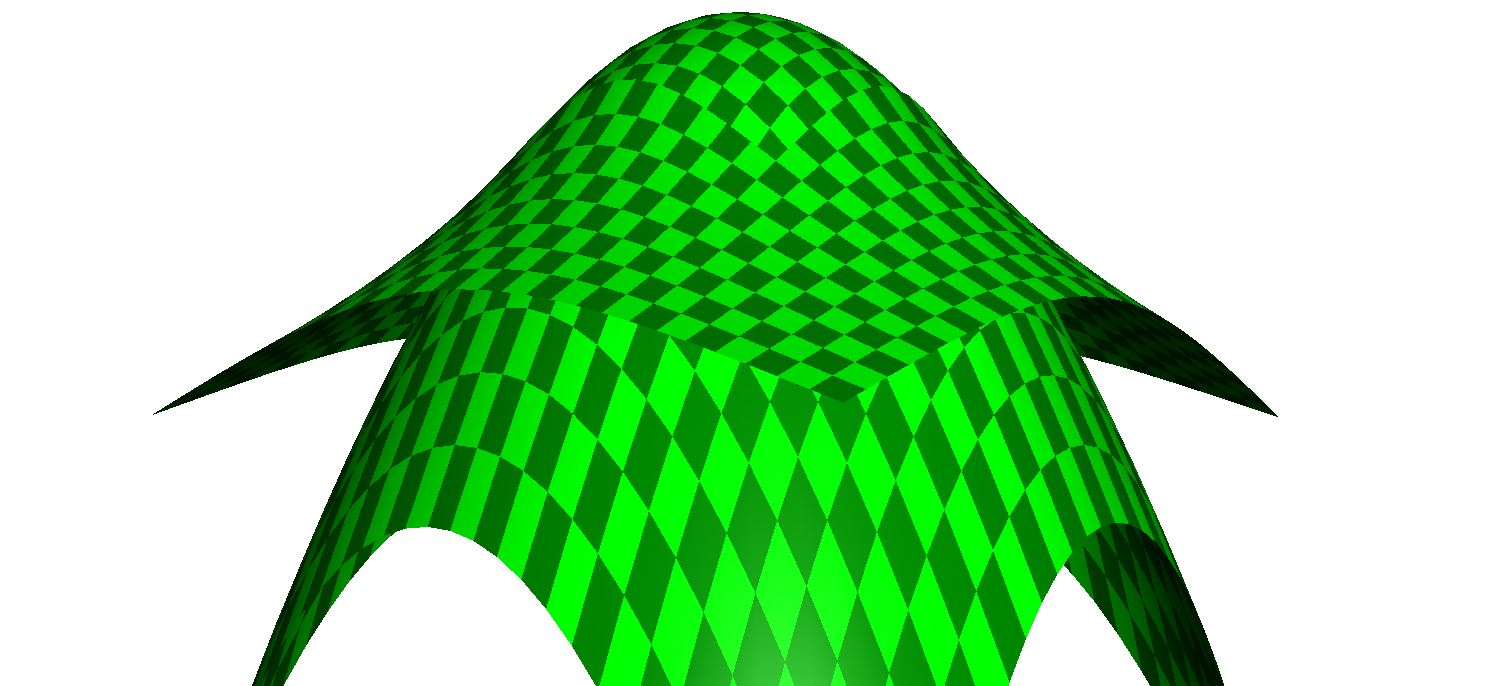}
  \end{center}
  \caption{Membrane problem in the presence of an obstacle.}
\end{wrapfigure}
While the existence and uniqueness of solution to the obstacle problem can be established exactly as in Problem \ref{Prob1}, its regularity theory is rather different. Indeed, solutions to Problem \ref{Prob1} are infinitely many differentiable, whereas the optimal regularity for the membrane restricted to lay above an obstacle drops to $C^{1,1}$; that is, the best one can hope is boundedness of second derivatives. This is indeed true and accounts an important Theorem firstly proven by J. Frehse (1972). 

The next, and rather more involved issue is to understand the smoothness of the interface between the contact set, $C:= \{u(x) = \varphi(x)\}$ and the non-contact set $\{u(x) > \varphi(x)\}$, the so called free boundary of the problem, $\Gamma$. Luis Caffarelli proved in \cite{Caff77} that if the contact set is thick enough in a neighborhood of a free boundary point $x_0$,   then $\Gamma$ is a $C^1$ surface around $x_0$. In slightly more precise terms, Caffarelli showed  the existence of a critical density $\varrho(r)$, such that if $x_0 \in \Gamma$ is such that for some $0< r \ll 1$, Width$\left (  {C} \cap B_r(x_0)\right ) \ge \varrho(r)$, then $\Gamma$ is a $C^{1,\alpha}$ (and thus $C^\infty$, by a result from \cite{KN}) surface around $x_0$. 

Four years later, in 1981, Caffarelli, partnering with H. Alt, publishes what would become a {\it magnum opus} of variational free boundary theory. In \cite{AC}, Alt and Caffarelli study regularity properties of non-negative local minimizers to the discontinuous functional, $J(u) := \int |\nabla u|^2 + \chi_{\{u> 0\}} dx \longrightarrow \text{min}.$ Here, $\chi_{\{u> 0\}}$ stands for the characteristic function of the set $\{u> 0\}$. Since this is a discontinuous function, one should expect the Laplacian of a minimizer to behave as a Dirac mass along the set of discontinuity, namely $\partial \{ u > 0 \} =: \Gamma$ --- the free boundary of the problem. This is indeed the case; Alt and Caffarelli show that a local minimizer,$ u$,  behaves linearly along the free boundary, in particular is Lipschitz continuous (the optimal regularity of the problem). They also show that, in some very week sense, the normal derivative of $u$ along the free boundary is constant. The most delicate part of the program is to show differentiability of the free boundary,  a $\mathcal{H}^{n-1}$-negligible set. This is accomplished by a rather elaborate implementation of methods along the lines {\it flatness implies regularity}, involving non-homogeneous blow-ups.

The investigation of sign changing minimizes of discontinuous functionals of Alt-Caffarelli type, as above, is motived by problems in the theory of jet flows, phase transmission, among others. From the mathematical perspective though, the analysis of  sign changing minimizes is rather more involved then of its one-phase counterpart. In particular, establishing Lipschitz estimates for such minimizers required a powerful new tool, namely a monotonicity formula, in the spirit of geometric measure theory. This is the contents of the celebrated work of Alt, Caffarelli and Friedman (1984), \cite{ACF}.

Within the theory of free boundary regularity, the motto {\it flatness implies regularity} attains the apogee in Caffarelli's trilogy \cite{CC1, CC2, CC3}, where he develops a rather complete existence and regularity theory for a very general class of  two-phase free boundary problems. Here, however, Lipschitz estimates for the free boundary presents itself as an intermediary step. Caffarelli's free boundary regularity slogan then becomes: {\it flatness implies Lipschitz, and Lipchitz implies differentiability}.  
 
\section{Back to diffusion}

Recently, a radical new geometric approach to the analysis of diffusive PDEs has been launched, in which degenerate points of ellipticity are seen as part of what has been termed ``non-physical free boundaries".  

Heuristically speaking, it is often in the realm of mathematics that the complexity of a given problem $\mathbb{P}$ is encoded within some special {\it entities} pertaining to it: singularities, bifurcations, degeneracies, blow-ups, discontinuity, sharp changes, etc. Significant advances on the problem $\mathbb{P}$ depend upon critical understanding of such distinct elements and how they affect the  order of the model. When it comes to the analysis of diffusive PDEs, geometric insights from the free boundary theory provide a rather powerful toolbox to investigate those special points; the so termed non-physical free boundaries.

As a way of example, the 1930's Schauder estimates bear the premise that the smoothness of the gradient (or the Hessian, depending on whether the problem is in divergence or non-divergence form) of a solution to a second order linear elliptic equation could never exceed the continuity of the medium. That is, if the coefficients, $\gamma(x)$, of a divergent form operator  $\ell \mapsto \div (\gamma(x) \ell )$ is $\alpha$-H\"older continuous, for some $0< \alpha < 1$, then the gradient of a solution to the homogeneous equation  $\div (\gamma(x) \nabla u) = 0$ is also  $\alpha$-H\"older continuous, for the same exponent $\alpha$. This is a celebrated result, which is far from being elementary, or even intuitive. Far less perceptive is the fact that the continuity of the gradient can be even superior to the continuity of the coefficients, $\gamma$ --- but only along its critical set $\mathscr{S}:=\{ \nabla u(x) = 0\}$, \cite{TT2, TT3}. 
\begin{wrapfigure}{r}{0.5\textwidth}
 \begin{center}
 \vspace{-0.4cm}
    \includegraphics[scale=0.25]{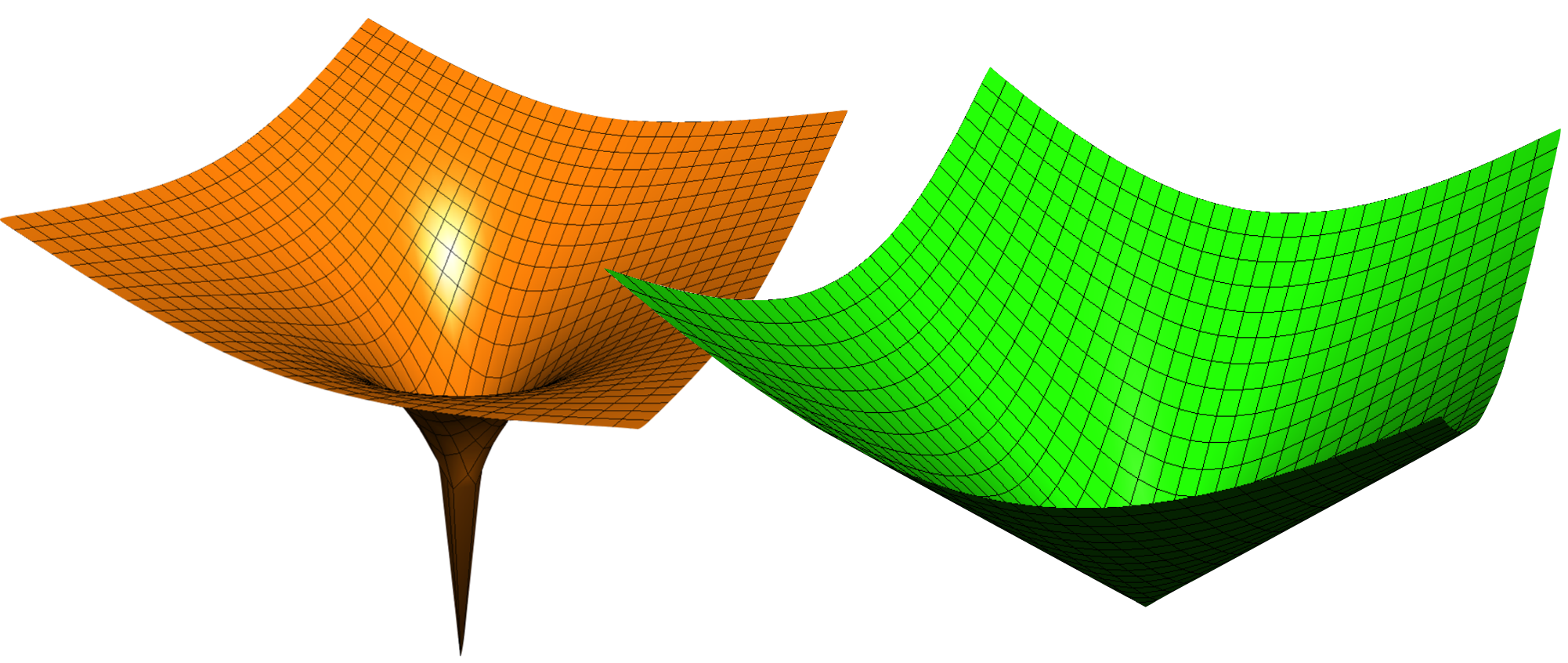}
  \end{center}
  \caption{Improved regularity along critical points: at a generic point, gradient (resp. hessian) of solutions may develop sharp cups; at a critical point, however, it can only appear much smoother corner-like singularities.}
\end{wrapfigure}

Such an improved regularity estimate becomes even more appealing in the context of degenerate elliptic equations, as in the theory of the $p$-Laplacian. This is because the critical set $\mathscr{S}:=\{ \nabla u(x) = 0\}$ is precisely the region in which the diffusion attributes of the operator collapse.  Striking enough, even if the medium does not have power oscillation decay, the gradient of a solution does, but only around points of $\mathscr{S}$.
 
The formal statements of such results are a bit too technical to be stated here; however it is noteworthy to comment that the core ideas for proving these theorems are genuinely geometric, and were largely influenced by the general free boundary framework mentioned above.  Several applications and enhancements of these methods have been successfully set forth in the past few years, leading to a plethora of other unanticipated results. This is currently a rather active line of investigation and it is likely that the analogy herein set will bear fruit in other branches of mathematical analysis.

\end{document}